\documentclass[11pt]{amsart}
\usepackage{amsmath}
\usepackage{amssymb}
\usepackage{amscd}

\theoremstyle{plain}
\newtheorem{thm}{Theorem}[section]
\newtheorem{lem}[thm]{Lemma}
\newtheorem{cor}[thm]{Corollary}
\newtheorem{prop}[thm]{Proposition}

\theoremstyle{definition}

\def \CPb {\overline{\mathbf{CP}}^{\,2}}

\def \PO {{\mathbf{CP}}^{\,1}}

\def \Z {\mathbb{Z}}
\def \Sig{\Sigma}

\def \SS {{S}^2\times {S}^2}

\def \vp {\varphi}

\def \a {\alpha}
\def \b {\beta}
\def \g {\gamma}
\def \d {\delta}
\def \k {\kappa}

\def \L {\Lambda}
\def \G {\Gamma}

\def \o {\omega}

\def \t {\tau}
\def \z {\zeta}

\def \bd {\partial}
\def \x {\times}
\def \- {\setminus}
\def \C {\subset}
\def \ve {\varepsilon}

\def \sw {\text{SW}}

\def \DD {\Delta}

\def \hZ {\hat{Z}}
\def \hY {\hat{Y}}
\def \ha {\hat{a}}
\def \sign {\text{sign}}

\def \bm {{\bf{m}}}

\begin{document}

\baselineskip.525cm
\title[Simply connected $4$-manifolds with the same   Seiberg-Witten
invariants]{Families of
simply connected $4$-manifolds with the same   Seiberg-Witten invariants}
\author[Ronald Fintushel]{Ronald Fintushel}
\address{Department of Mathematics, Michigan State University \newline
\hspace*{.375in}East Lansing, Michigan 48824}
\email{\rm{ronfint@math.msu.edu}}
\thanks{The first author was partially supported NSF Grant DMS0072212
and the second author by
NSF Grant DMS9971667}
\author[Ronald J. Stern]{Ronald J. Stern}
\address{Department of Mathematics, University of California \newline
\hspace*{.375in}Irvine,  California 92697}
\email{\rm{rstern@math.uci.edu}}
\date{October 14, 2002}

\maketitle

\section{Introduction\label{Intro}}

In 1994, the introduction of the Seiberg-Witten equations quickly
fostered optimism in
4-manifold topology. Invariants derived from these equations
immediately led to the solutions
of several outstanding conjectures, and it was felt that  the
classification of 4-manifolds
was finally in sight. After eight years, however, the opposite seems to be
true. Topological
constructions, along with  the Seiberg-Witten invariants, have
demonstrated that simply
connected smooth 4-manifolds are much more complicated than earlier
envisioned, and the
classification of smooth $4$-manifolds has retreated beyond the
visible horizon.  In this
paper we add to this quagmire and exhibit constructions  which yield
new infinite families of
homeomorphic simply connected 4-manifolds, all  of which have the
same Seiberg-Witten
invariants. These manifolds have many  characteristics which lead to
the belief that they
ought not be diffeomorphic, but this remains unsettled. We also
introduce a construction which
readily converts an irreducible symplectic $4$-manifold $X$ to a
family of infinitely many
irreducible nondiffeomorphic but mutually homeomorphic smooth $4$-manifolds. Exactly one member of this family will admit a symplectic structure.

There are many examples of the first sort
which are already known.
For  example, for each odd integer $r\ge3$, one can consider the
Horikawa surfaces with
holomorphic euler number $2r-1$ and $c_1^2=4r-8$. (See e.g.
\cite{GS}.) For each such
$r$, there are two deformation classes of these simply connected
complex algebraic  surfaces,
and the two deformation types are homeomorphic. Since the Horikawa
surfaces are general type,
their Seiberg-Witten invariants are $\sw = t_K - t^{-1}_K$  where $K$
is the canonical class
(which is primitive). (We shall view the  Seiberg-Witten invariants
of simply connected
$4$-manifolds as elements of the integral group ring
$\Z H_2(X;\Z)$. The notation $t_b$ will be used to denote the element
of the group ring
corresponding to $b\in H_2(X;\Z)$ and its (integer) coefficient will
at times be denoted by
$\sw_X(b)$.) Thus the Seiberg-Witten invariants fail to  distinguish
these manifolds. Whether
or not they are indeed diffeomorphic is an extremely interesting open question.

Knot surgery \cite{KL4M,ICM} can be used to produce infinite families
of homeomorphic simply
connected 4-manifolds with the same  Seiberg-Witten invariants. For
example, if we produce a
family of knots $\{ K_n\}$ all of which  have the same symmetrized
Alexander polynomial
$\DD_{K_n}(t)$ (for example, take an  infinite list of knots of
Alexander polynomial equal to
$1$) then the manifolds
$X_n$ obtained from knot surgery on an elliptic fiber of the
K3-surface will all be
homeomorphic to K3 and will have the same Seiberg-Witten invariants.
Although these examples
seem likely to be nondiffeomorphic, one needs to use care here. S.
Akbulut \cite {Ak} has
shown that for any knot $K$ that the knot surgery manifolds  arising
from K3 by using $K$ and
its mirror image
$-K$ are diffeomorphic. The proof uses special symmetry properties,
but it indicates that one
should not assume too quickly that the knot  surgery manifolds
arising from different knots
ought to be nondiffeomorphic.

In this paper we shall produce three further families. The first is a
generalization of the
authors' construction \cite{1bc} of nonsymplectic manifolds with one
basic class. For this
generalization we begin by constructing families of homeomorphic
simply connected symplectic
manifolds with one basic class with Seiberg-Witten invariant $\pm 1$,
i.e. with $\sw = t_K \pm
t^{-1}_K$. As evidence that these manifolds with one basic class are
never diffeomorphic to a complex surface, we show that this
construction yields both families of manifolds homeomorphic to a
simply connected complex surface as well  as families of manifolds
that are homeomorphic to no complex surface.  

We then twist this construction to yield further families whose Seiberg-Witten invariants 
have similar properties. In addition, in each manifold of this family, we find nullhomologous tori with the property
that a $(-1/m)$-logarithmic transform on any of them multiplies the
Seiberg-Witten invariant by $m+1$.
This results in families of homeomorphic simply
connected 4-manifolds with $\sw = (m+1)(t_K \pm t^{-1}_K)$. We
conjecture that no two of these
manifolds are diffeomorphic. Note that if $m \ne 0$ these
manifolds are nonsymplectic.
Thus, this construction shows how one can perform surgery to
make a symplectic manifold nonsymplectic.

The third collection of families builds on a  construction of the
authors which gives a
counterexample to the existence of a symplectic Parshin-Arakelov
theorem. In particular, we
show that for a fixed even (resp. odd) integer $g$ there are
infinitely many distinct genus
$g$ Lefschetz fibrations on a homotopy rationally elliptic (resp. K3)
surface. These examples
have been promised for some time and were first presented in a  talk
of the first author at
the 1998 Aarhus Topology Conference.

\section{Construction I}

Fix a simply connected symplectic $4$-manifold $X$, and let $C$ be a
smoothly embedded
symplectic surface in $X$ which has genus $n\ge 2$ and  self-intersection
$0$. Given $X$, a
positive integer $g$, and a particular  genus $g$ fiber bundle $Y$
over a genus $n$ surface
constructed below, we shall associate a symplectic simply connected
$4$-manifold $Z=Z(X,C,g)$
with $c_1^2(Z)=c_1^2(X)+8g(n-1)$, $\chi(Z)=\chi(X) + g(n-1)$, and
with $\sw = t_K - t^{-1}_K$.
Here $\chi = \chi(X)$ is
one-fourth the sum of the Euler
characteristic and the signature of $X$, which is the is the
holomorphic Euler characteristic
in the case that $X$ is a complex
surface.
We shall choose a family of
manifolds $X_j$,
homeomorphic to $X$, and distinguished by their Seiberg-Witten
invariants.  Each $X_j$
contains a genus $n$ symplectic surface $C_j$ with self-intersection
$0$.  We fiber sum along
these surfaces to obtain the family $Z(X_j,C_j,g)$ with the same
Seiberg-Witten invariants. 

To construct the
bundle $Y$, let $K_g$ denote the $(2g+1, -2)$-torus knot, pictured in
Figure~1, and let
$M_{K_g}$ denote the $3$-manifold obtained by performing $0$-framed surgery on
$K_g$. This manifold has the integral homology of $S^2\x S^1$. In
Figure 1 we can see an
obvious genus $g$ Seifert surface for $K_g$. In $M_{K_g}$, fix  a closed genus
$g$ surface $\Sig_g$ obtained from capping off this Seifert surface
with a disk.

\centerline{\unitlength .75cm
\begin{picture}(6,9)
\put (2,4.5){\oval(3,7)[l]}
\put (3,4.5){\oval(3,7)[r]}
\put (2,8){\line(3,-4){1}}
\put (3,6.66){\line(-3,-4){.4}}
\put (3,8){\line(-3,-4){.4}}
\put (2,6.66){\line(3,4){.4}}
\put (2,6.66){\line(3,-4){1}}
\put (2,5.33){\line(3,4){.4}}
\put (2.33,3){$\bullet$}
\put (2.33,3.67){$\bullet$}
\put (2.33,4.33){$\bullet$}
\put (2,2.33){\line(3,-4){1}}
\put (3,2.33){\line(-3,-4){.4}}
\put (2,1){\line(3,4){.4}}
\put (1.6,.35){Figure 1}
\end{picture}}

Now fix a $g\ge 1$, and set $K=K_g$ and $\Sig = \Sig_g$. Because $K$
is a fibered knot,
$S^1\x M_K$ is a symplectic manifold which is fibered over $T^2$ with
symplectic fiber
$\Sig$. The fiber sum
\[ Y_{2,g}= S^1\x M_K\,\#_{\Sig}\, S^1\x M_K \]  of two copies of
$S^1\x M_K$ is again a symplectic manifold.  This manifold has euler number
$e(Y_{2,g})=4g-4$ and signature
$\sign(Y_{2,g}) = 0$. In the $j$th copy of $S^1\x M_K$, let
$T_j$ denote the torus $S^1\x m$ where $m$ is a meridian to $K$ (and
thus generates
$H_1(M_K;\Z)$). This torus is a section to the fibration of
$S^1\x M_K$. It follows from work of Meng and Taubes
\cite{MT} that any basic class of $S^1\x M_K$ must be a multiple of
$T_j$. They show that the
Seiberg-Witten invariant of $S^1\x M_K$ is
$\DD_K(t^2)/(t-t^{-1})^2$ where
$t=t_{T_j}\in \Z H_2(S^1\x M_K;\Z)$ and $\DD_K$ is the symmetrized
Alexander polynomial of
$K$. Note that the intersection number of
$T_j$ and $\Sig$ is $T_j\cdot\Sig=1$.

We wish to determine the Seiberg-Witten invariant of $Y_{2,g}$. Write
\[ Y_{2,g} = ((S^1\x M_K)\- (\Sig\x D^2))\cup ((S^1\x M_K)\- (\Sig\x
D^2)).\] In
$S^1\x M_K$, the torus $T_j$ intersects $\Sig\x D^2$ in some
$\{y_j\}\x D^2$. We may assume that the identification of the two
copies  of $\Sig$ in the
fiber sum is chosen so that $y_1$ is identified with $y_2$. We  then
obtain a genus 2 surface
$(T_1\- D^2)\cup  (T_2\- D^2)$ in the fiber sum. (Of course,
$Y_{2,g}$ is a fiber bundle over a genus 2 surface, and we have just
constructed a section.) Let
$\t_{1,2}$ denote the class of this surface in
$H_2(Y_{2,g};\Z)$. Notice that $\t_{1,2}\cdot\Sig=1$ (and that
$\t_{1,2}^2=0=\Sig^2$).

There are also $4g$ classes in $H_2(Y_{2,g};\Z)$ which are obtained as
follows:  For a fiber
$\Sig\x \{t\}\C \bd(\Sig\x D^2)$, the inclusion
$H_1(\Sig\x \{t\};\Z)\to H_1((S^1\x M_K)\- (\Sig\x D^2);\Z)$ is
trivial.  Thus, the
identification of fibers and a collection of loops $\{ a_i\}$ on
$\Sig$ which gives a basis for $H_1(\Sig;\Z)$
 gives rise to $2g$  2-dimensional homology classes $V_i$ formed
from unions of bounding
surfaces.  We shall refer to these as {\bf{vanishing classes}}. For each $a_i$ there is a rim torus $R_i=  a_i\x\bd D^2$.
With appropriate
orientation choices, we have for all
$i,j$:
\begin{enumerate}
\item[] $R_{j}\cdot V_{i} = a_i\cdot a_j$ \hspace{.15in} (the
skew-symmetric intersection form
on $H_1(\Sig;\Z)$),
\item[] $R_j\cdot \Sig=0$
\item[] $R_j\cdot \t_{1,2} =0$
\end{enumerate} The intersection matrix $Q$ of $\{ V_{i},R_{j}|
i,j=1,\dots,2g\}$ has
determinant 1, and these classes along with $\t_{1,2}$ and $\Sig$
form a basis for
$H_2(Y_{2,g};\Z)$.

\centerline{\unitlength 1cm
\begin{picture}(12,6)
\put (2,3){\oval(3,4)[l]}
\put (3,3){\oval(3,4)[r]}
\put (2,5){\line(3,-4){1}}
\put (3,3.66){\line(-3,-4){.4}}
\put (3,5){\line(-3,-4){.4}}
\put (2,3.66){\line(3,4){.4}}
\put (2,3.66){\line(3,-4){1}}
\put (2,2.33){\line(3,4){.4}}
\put (2,2.33){\line(3,-4){1}}
\put (3,2.33){\line(-3,-4){.4}}
\put (2,1){\line(3,4){.4}}
\put (2.4,3.66){\oval(1.5,1.4)[l]}
\put (2.6,3.66){\oval(1.5,1.4)[r]}
\put (1.5,4.25){\small{$a_1$}}
\put (1.85,.5){Figure 2}
\put (9,3){\oval(3,4)[bl]}
\put (7.49,3){\line(0,1){.125}}
\put (10,3){\oval(3,4)[r]}
\put (9,5){\line(-1,0){.85}}
\qbezier (7.49,4.35)(7.6,4.9)(8.15,5)
\put (9,5){\line(3,-4){1}}
\put (10,3.66){\line(-3,-4){.4}}
\put (10,5){\line(-3,-4){.4}}
\put (9,3.66){\line(3,4){.4}}
\put (9,3.66){\line(3,-4){.265}}
\put (9,2.33){\line(3,4){.4}}
\put (9,2.33){\line(3,-4){1}}
\put (10,2.33){\line(-3,-4){.4}}
\put (10,2.33){\line(-3,4){.6}}
\put (9,1){\line(3,4){.4}}
\put (7.1,3.2){\line(1,0){2.3}}
\put (7.65,4){\line(1,0){1.5}}
\put (7.1,4){\line(1,0){.25}}
\qbezier (9.4,3.2)(10,3.6)(9.4,4)
\qbezier (7.1,3.2)(6.5,3.6)(7.1,4)
\put (7.49,4.35){\line(0,-1){1}}
\put (6.4,3.75){\small{$a_1^+$}}
\put (8.85,.5){Figure 3}
\end{picture}}

It will be useful to have a more specific description of the vanishing 
classes $V_j$.  As we have mentioned, the complement in $S^3$ of the
$(2g+1, -2)$ - torus knot $K_g$, fibers over the circle, and its
fiber is a punctured genus $g$ surface, $\Sig_g'$. If we make the
obvious choices for the loops $\{ a_j\}$ which are mentioned above,
their associated Seifert linking pairing is
\[ \ell k(a_i,a_j^+)=
\begin{cases}1,\ \ & j=i, i+1 \\0, &
{\rm{otherwise}}\end{cases}
\]
where $a_j^+$ is the positive push-off
of $a_j$. Notice that $a_j^+$ bounds a genus 1 surface  in $S^3\-
K_g$.  A typical element $a_1$ of this basis and its push-off $a_1^+$
for
$K_1$, the left-hand trefoil knot, are shown in Figures 2 and 3.

Again writing $K$ for $K_g$, we consider the fiber bundle
\[ p:
(S^1\x M_K)\- (\Sig\x D^2)\to T^2\- D^2. \]
The base deformation
retracts to a wedge of two circles, $S^1_1\vee S^1_2$, where
$p^{-1}(S^1_1)=
S^1_1\x \Sig$  and $p^{-1}(S^1_2)=M_K$. Let
$S^1_1\cap S^1_2 = \{x_0\}$, and let $y_0$ be a point in the boundary
of $T^2\- D^2$. Let $\d: [0,1]\to T^2\- D^2$ parametrize $S^1_2$ so
that $\d(0)=x_0$ and so that the positive push-off $a_i^+$ lies in
$p^{-1}(\d(.1))$.

Now let $\g$ be an embedded path in $T^2\- D^2$
with $\g(0)=y_0$ and  $\g(1)=\d(.1)$. In $p^{-1}(y_0)\cong\Sig$
consider the circle $a_{i,\bd}$, the obvious circle in the fiber over
$y_0$. This is a boundary component of an embedded annulus in
$p^{-1}(\g)$ whose other boundary component is  $a_i^+$, which bounds
a genus 1 surface in
$M_K=p^{-1}(S^1_2)$. We thus obtain in $(S^1\x
M_K)\- (\Sig\x D^2)$ a genus 1 surface $G_{i,1}$ with boundary $a_i$.
Similarly, in the second copy of $(S^1\x M_K)\- (\Sig\x D^2)$, we
have a similar genus 1 surface $G_{i,2}$.

So far, we have left ourselves
the freedom of choosing the gluing map
which constructs the fiber sum. At this point, we specify that the
chosen gluing must take the boundary of $G_{i,1}$ diffeomorphically
onto the boundary of $G_{i,2}$ for each $i$. Then $V_i= G_{i,1}\cup
G_{i,2}$ is a genus 2 surface in $Y_{2,g}$.

\begin{lem} The genus 2
surface $V_i$ has self-intersection number $2$ and satisfies
$V_i\cdot\Sig=0$, $A_i \cdot\t_{1,2}=0$, and
$V_i\cdot R_{j}= a_i\cdot a_j$. \end{lem}
\begin{proof} Except for
the calculation of the self-intersection number, everything is clear
from the construction of $V_i$. Let ${S^1_2}'$ denote a circle in
$T^2\- D^2$ which is disjoint from and parallel to $S^1_2$, and let
$\g'$ be a path parallel to $\g$ from a point in $\bd (T^2\- D^2)$ to
a point on ${S^1_2}'$. Using the structure of the fibration,
$G_{i,1}$ in        $p^{-1}(\g\cup S^1_2)$ can be pushed onto
$G_{i,1}'$ in $p^{-1}(\g'\cup {S^1_2}')$.
Depending on the choices
made for ${S^1_2}'$ and $\g'$, either $S_1^2$ will intersect $g'$ in
a single point or ${S^1_2}'$ will intersect $\g$ in a single positive
point. In either case the only fiber in which $G_{i,1}$ and
$G_{i,1}'$ intersect is the fiber over this point. The intersection
of each surface with the fiber is a push-off of $a_i$, and these two
push-offs intersect once on this fiber.
The same is true for surfaces
$G_{i,2}$ and $G_{i,2}'$, implying that $V_i^2=2$.
\end{proof}

The adjunction inequality implies that any basic class must intersect
a smoothly embedded
surface of genus $h$ and self-intersection $2h-2$ with intersection
number equal to $0$. In
particular, the intersection number of any  basic class of
$Y_{2,g}$ with a rim torus must be $0$. Consider a basic class $k$ of
$Y_{2,g}$ and write
\[ k=t\,\t_{1,2}+s\,\Sig+\sum_{i=1}^{2g}u_i\,R_i+v_i\,V_i \] The intersections
$k\cdot R_j=0$, $j=1,\dots,2g$, give rise to the equation
$Q^T{\bf{v}}={\bf{0}}$ where $Q^T$ is the transpose of the
intersection matrix $Q$ and
${\bf{v}}=(v_1,\dots,v_{2g})$. Since $Q$ is nonsingular, all $v_i=0$.

We see that a basic class of $Y_{2,g}$ has the form
\[ k=t\,\t_{1,2}+s\,\Sig+\sum_{j=1}^{2g}u_j\,R_j \] Another application of
the adjunction  inequality gives $0=V_i\cdot k$, and we get the equation
$Q{\bf{u}}={\bf{0}}$. Once again this means that each coefficient
$u_j=0$. Thus each basic
class of $Y_{2,g}$ has the form $k=t\,\t_{1,2}+s\,\Sig$.

Since $Y_{2,g}$ is a symplectic manifold, it has simple type; that is, for each
basic class $k$,
\[ k^2 = 3\,\sign (Y_{2,g}) +2 e(Y_{2,g}) = 8g-8\] So if
$k=t\,\t_{1,2}+s\,\Sig$ then
$2st=8g-8$. Applying the adjunction inequality to $\t_{1,2}$ and to
$\Sig$ yields:
\[2\ge \t_{1,2}^2+|k\cdot\t_{1,2}|=|s| \hspace{.15in} {\text{and}}
\hspace{.15in} 2g-2\ge \Sig^2+|k\cdot\Sig |=|t|\] Thus if $g>1$, the
only basic classes  of $Y_{2,g}$
are
$\pm \b$ where $\b= (2g-2)\,\t_{1,2}+2\,\Sig$. Clearly, $\b$ is the
canonical class of $Y_{2,g}$, and
\[{\sw}_{Y_{2,g}}=t_{\b}+(-1)^{g-1}t_{\b}^{-1}\]
In the special case
$g=1$, $Y_{2,1}$ is a torus bundle over a genus 2 surface. In this
case the above argument only shows that for the canonical class 
$\b=2\Sig$ of $Y_{2,1}$,
\[{\sw}_{Y_{2,1}}=t_{\b}+c+t_{\b}^{-1}\]
where $c=\sw_{Y_{2,1}}(0)$.

In general, let $Y_{n,g}= (S^1\x M_K)\#_{\Sig}(S^1\x M_K)\#_{\Sig}\cdots
\#_{\Sig}(S^1\x M_K)$ \  ($n$ copies) where $K=K_g$.  Then $Y_{n,g}$
is a
symplectic manifold which is a fiber bundle
with a genus $g$
fiber and a genus $n$ base. It contains a section $S$
of square $0$
obtained from gluing together the punctured copies of
$T_i=S^1\x m_i$. We still denote the genus $g$ fiber by $\Sig$.
Thus
$S\cdot\Sig=1$. The same analysis as above  shows that for $g>1$ the
only basic
classes of $Y_{n,g}$ are
$\pm\b$ where $\b=(2g-2)S+(2n-2)\Sig$ is the canonical class of $Y_{n,g}$, and
\[ \sw_{Y_{n,g}}= t_{\b}+(-1)^{(g-1)(n-1)}t_{\b}^{-1} \]
Also
$\sw_{Y_{n,1}}$ is a monic symmetric Laurent polynomial in powers of
$t_{\Sig}$, and the highest power that occurs is $\pm (2n-2)$. These
terms correspond to $\pm \b$ where $\b=(2n-2)\Sig$ is the canonical
class.

We note for use
below that when $n>2$, we have $Y_{n,g}=\hY_1\cup\dots\cup \hY_n$, where each
$\hY_i$, $i\ne 1,n$, is a bundle over $T^2\- (D^2\cup D^2)$.
In these
latter manifolds there are punctured tori $G_{i,j}^{\pm}$ so that we
get genus 2, self-intersection 2 (vanishing) surfaces $V_{i,j}=G_{i,j}^+\cup
G_{i,j+1}^-$, for $i=1,\dots,2g$ and $j=1,\dots, n-1$.

Recall that the data given for our construction consists of a simply
connected symplectic
$4$-manifold $X$ with an embedded symplectic surface $C$ of genus
$n\ge 2$ and self-intersection number
$0$, and the integer $g\ge 1$. Form the  symplectic manifold $Z(X,C,g)$ as
the fiber sum of $X$ and $Y_{n,g}$:
\[ Z=Z(X,C,g)= X\#_{C=S}Y_{n,g}.\]

\begin{prop} \label{p} If $\pi_1(X\- C)=1$, then $\pi_1(Z)=1$. \end{prop}
\begin{proof} The $i$th $\pi_1(S^1\x M_K)$ is normally generated by
the image of
$\pi_1(T_i)$. Since $T_i$ intersects $\Sig$ in a single point, a
normal circle to
$\Sig$ lies in $T_i$. Thus
$\pi_1(S^1\x M_K\- \Sig)$ is normally generated by the image of
$\pi_1(T_i\- {\text{point}})$. An inductive application of Van
Kampen's theorem shows that
$\pi_1(Y_{n,g})$ is normally generated by the image of $\pi_1(S)$. Thus if
$\pi_1(X\- C)=1$, we have
$\pi_1(Z)= \pi_1(Y_{n,g}\- S)/\pi_1(S\x S^1) =1$.
\end{proof}

\begin{lem}\label{rim} Let $S$ be an orientable surface in the
compact orientable 4-manifold $Y$. The subgroup of $H_2(Y\- S;\Z)$ generated
by the rim tori of $S$ is isomorphic to $H^1(S)/{\rm{im}}(H^1(Y))$.
It follows that the rim tori of $Z(X,C,g)$ arising from loops in
$C=S$ are all nullhomologous in
$Z(X,C,g)$. \end{lem}
\begin{proof} The group generated by the classes of  rim tori is the
kernel of    $H_2(Y\- S)\to H_2(Y)$. Let $N$ be a tubular
neighborhood of $S$. Then the main statement of the lemma follows
from
\[ \begin{CD}  H_3(Y)\to H_3(Y,Y\- S) @>{\bd}>> H_2(Y\-
S)\to H_2(Y) \end{CD} \]
once we note that $H_3(Y,Y\- S)\cong H_3(N,
\bd N)$ and invoke Poincar\'e duality. The second statement follows
since $H^1(Y_{n,g})\cong H^1(S)$.
\end{proof}

Dually, in the situation of the lemma, there is a 
homomorphism $H_3(Y)\to H_1(S)$ given by intersection. The lemma 
states that the group generated by the rim tori is the quotient of 
$H_1(S)$ by the image of $H_3(Y)$. Note that in the special case 
where $p:Y\to B$ is a surface bundle over a surface, and $S$ is the 
image of a section, this implies that the subgroup of $H_2(Y\- S)$ 
generated by the rim tori of $S$ vanishes.

For a basic class $\k$ of $X$, the adjunction inequality implies that
the maximal
intersection of $\k$ with $C$ is $\k\cdot C = 2n-2$. This
intersection number is achieved
when $\k$ is the canonical class of $Z$, and  whenever this maximal
intersection number is
achieved, there is a surface $B_{\k}$ in $X$,  representing
$\k$, which intersects $C$ positively in exactly $2n-2$ points. We
can then form the class
$\b_{\k}\in H_2(Z,\Z)$ which is represented by the union of $B_{\k}$
with $2n-2$ disks
removed and a smooth surface representing $\b$ with the $2n-2$
normal disks removed at the
points where it intersects $S$. Any basic class $\a$ of $Z$ satisfies
\[\a^2= c_1^2(Z)=c_1^2(X)+c_1^2(Y_{n,g})+8n-8\] Note that
$\ve_{\k}=\b_{\k}+2\,C$
has exactly this square.

It follows from \cite{MST} (and the fact that there are no rim tori) that
\[ \sw_Z(\ve_{\k})=\sw_{Y_{n,g}}(\b)\cdot\sw_X(\k)= \sw_X(\k) \] and that
these are the only basic
classes of $Z$ which have intersection number $2n-2$ with $C$.
To express this, write the Seiberg-Witten  invariant of $Z$ as
      \[\sw_Z=\sw_{\{Z;C,-\}}+\sw_{\{Z;C,{\text{max}}\}}\]  where the the
first summand consists
of terms of the form $b\,t_{\ell}$ where $|\ell\cdot C|<  2n-2$, and
the latter summand
consists of terms corresponding to basic classes whose intersection
number with $C$ is
$\pm (2n-2)$. Similarly write
\[\sw_X=\sw_{\{X;C,-\}}+\sw_{\{X;C,{\text{max}}\}}.\] Then we have:

\begin{prop}\label{Z} Let $X$ be a symplectic 4-manifold containing
an embedded symplectic
surface $C$ of genus $n$ and self-intersection $0$.  Suppose also that
$\pi_1(X\- C)=0$. Then for each $g\ge 1$,
$Z=Z(X,C,g)$ is simply connected, and if the Seiberg-Witten invariant of $X$ is
\[\sw_X=\sw_{\{X;C,-\}}+\sum_{k\cdot C=2n-2} c_k\,(t_{k}+
(-1)^{\chi(X)}t_{k}^{-1})\]
then the the Seiberg-Witten invariant of
$Z$ is
\[\sw_Z=\sw_{\{Z;C,-\}}+\sum_{k\cdot C=2n-2}c_k\,(t_{\ve_{k}}+
(-1)^{\chi(Z)}t_{\ve_{k}}^{-1})\]
\end{prop}

It seems quite likely that generally one has $\sw_{\{Z;C,-\}}=0$.
Under special hypotheses on
the manifold $X$, one can verify this. For example, let
$X$ be the the knot surgery manifold obtained from the K3-surface by
replacing a  neighborhood
of a torus fiber $F$ with the manifold $S^1\x(S^3\- K')$, where
$K'$ is a fibered knot of genus $n-1$. In the K3-surface let $C'$ be
a symplectically embedded
torus homologous to a fiber plus a section.  Then in $X$, a disk in
$C'$ is replaced with a
Seifert surface for $K'$, and we get a  symplectic surface $C$ of genus
$n$ and self-intersection $0$. In
\cite{1bc} it is shown that $X\- C$ is simply connected, and as in
that paper, adjunction
inequality arguments can be used to show that
$\sw_{\{Z;C,-\}}=0$, and $\sw_Z=t_F^{n-1}+t_F^{1-n}$.
These constructions give the examples promised at the beginning of this section.

\section {Other interesting $Z(X,C,g)$}

In this section  we shall choose $X$ to be a complex surface with $C$
a holomorphically
embedded submanifold. For a finite number of $g$ we will show that
$Z(X,C,g)$ has $\sw = t_K -
t^{-1}_K$, and by a theorem of Persson, Peters, and Xiao
\cite{complex}, its homeomorphism
type does not support any complex structure. For the remaining $g$
there are indeed complex
manifolds homeomorphic to $Z(X,C,g)$. We conjecture that none of the
$Z(X,C,g)$ have a
compatible complex structure.  The relevant result is the following
somewhat surprising restriction on the geography of complex surfaces
which support a spin
structure, i.e. with vanishing second Stiefel-Whitney class.

\begin{thm} [Persson, Peters, and Xiao \cite{complex}]
\label{restrict} Let $X$ be a simply
connected spin surface whose characteristic numbers satisfy
\[ 2\chi \le c_1^2 < 3(\chi-5)\] then $c_1^2=2(\chi -3)$ with $c_1^2 =
8k$ and $k$ odd or
$c_1^2=\frac83(\chi-4)$ with $\chi =1$ (mod 3).
\end{thm}

Let $H(m)$ be the spin Horikawa surface with $\chi=8m -1$ (and
$c_1^2(H(m))=2\chi-6$). It is
well-known that each $H(m)$ supports a genus $2$ Lefschetz fibration.
Let $C$ be a fiber.
Then  the manifolds $Z(m,g)=Z(H(m),C,g)$ are spin and

$$c_1^2(Z(m,g))=16m+8g-8$$
$$ \chi (Z(m,g))=8m+g-1$$

Thus, whenever we fix $m$ and choose $g$ so that the characteristic
numbers of $Z(m,g)$ are
restricted by  Theorem~\ref{restrict}, we obtain symplectic manifolds
that support no complex
structure. For the remaining $g$, there are complex manifolds in the
homeomorphism type of
$Z(m,g)$; however, we conjecture that no $Z(m,g)$ supports a complex structure.

One can verify that $X(m,g)$ supports no complex structure whenever
$g<8/5m-2$ and if $m=1$
(mod 3) then $g \ne 0$ (mod 3), if $m=2$ (mod 3) then $g \ne 1$ (mod
3), or if $m=0$ (mod 3)
then $g \ne 2$ (mod 3).

For example, there are no restricted examples starting with $H(1)$
and $H(2)$, i.e. $Z(1,g)$
and $Z(2,g)$ are homeomorphic to complex surfaces. However $Z(3,g)$
is restricted when $g=1$,
$Z(4,g)$ is restricted when $g=1,2,4$, $Z(5,g)$ is restricted when
$g=2,3$, Z(6,g) is
restricted when $g=1,3,4,6,7$, etc.

Note that since $H(m)$ is the $m$-fold fiber sum of $H(1)$ along a
genus $2$ fiber, $Z(m,g)$
is the fiber sum of $Z(1,g)$ with $m-1$ copies of
$H(1)$.

\section {Construction II}

This construction is similar to the last, but has properties which
will be useful when performing the surgeries described in the next
section. Again, let us begin with a
simply connected symplectic 4-manifold $X$ containing an embedded
symplectic surface $C$ of genus $n\ge 2$ and self-intersection $0$. For
any $g\ge 1$ we can  form the manifold $Z(X,C,g)$ as above. Suppose
that we tweak this construction:
As a warm-up, we begin with a simple case. Let $K=K_g$ and
form a 
twisted fiber sum as follows:
\[ Y'_{1,g,n-1} = (S^1\x 
M_K)\#_{\Sig_g=S_g}Y_{g,n-1} =Y_{1,g}\#_{\Sig_g=S_g}Y_{g,n-1}.\]
The 
genus $g$ section $S_g$ of $Y_{g,n-1}$ is identified with the genus 
$g$ fiber of
$S^1\x M_K$ in the fiber sum. This new symplectic 
manifold
contains the genus $n$ surface $S'$ of self-intersection $0$ obtained
from the sum of the genus $n-1$ fiber of $Y_{g,n-1}$  together with 
the genus 1 section $S_1$ of $Y_{1,g}$. Define $Z'_{1,g,n-1}(X,C)= 
X\#_{C=S'}Y'_{1,g,n-1}$.

\begin{prop} \label{pi1Z'} If $X\- C$ is simply connected, then so is $Z'_{1,g,n-1}(X,C)$.
\end{prop}
\begin{proof} As in the proof of Proposition~\ref{p},
$\pi_1(Y_{1,g}\- \Sig_g)$ is
normally generated by the image of $\pi_1(S_1\- {\text{pt}})$. Also, 
$\pi_1(Y_{g,n-1}\- S_g)$ is normally generated by the image of the
fundamental group of a section $S_g'$, disjoint from $S_g$, and by
the normal circle to $S_g$ (which lies on a fiber).
In $\pi_1(Y'_{1,g,n-1})$, the image of
$\pi_1(S_g')$ is identified with $\pi_1(\Sig_g)$ in
$\pi_1(Y_{1,g}\- \Sig_g)$; thus $\pi_1(Y'_{1,g,n-1})$ is normally
generated by the image of $\pi_1(S')$. As in Proposition~\ref{p},
if $\pi_1(X\- C)$ vanishes, it
follows that $Z'_{1,g,n-1}(X,C)$ is simply connected.
\end{proof}

Note that
$H_2(Y_{g,n};\Z)$ has rank $4g(n-1)+2$, whereas the rank of
$H_2(Y'_{1,g,n-1};\Z)$ is
$4(g-1)(n-1)+2$. The point is that the fiber sum along $\Sig_g=S_g$
contributes no rim tori classes (nor the associated vanishing classes) because rim tori to the
section $S_g$ bound in $Y_{g,n-1}\- S_g$. (See the remark following 
Lemma~\ref{rim}.)

The fiber $\Sig_{n-1}$ of $Y_{g,n-1}$ has a basis for $H_1$ represented by the loops 
$a_1,\dots, a_{2n-2}$ which were discussed in \S2.  Using the inclusion $\Sig_{n-1}\C S'$ and the identification of $S'$ with $C$ in $Z'_{1,g,n-1}(X,C)$ we obtain loops $\bar{a}_i$ on $C\x {\text{\{point\} }}$ in $C\x \bd D^2\C X\- (C\x D^2)$.
In the fiber sum with $X$, the rim tori to the $\bar{a}_i$ and the vanishing classes contribute $2n-2$ new hyperbolic pairs. Thus $Z(X,C,g)$ and $Z'_{1,g,n-1}(X,C)$ have isomorphic
intersection forms, and we get:

\begin{prop} The symplectic manifolds
$Z(X,C,g)$ and $Z'_{1,g,n-1}(X,C)$ are homeomorphic.
\end{prop}

Let $\b'= (2n-2)\Sig_g+(2g-2)S'$ be the canonical class of $Y'_{1,g,n-1}$. An argument
similar to that of \S2 shows that
\[ \sw_{Y'_{1,g,n-1}}=t_{\b'}+t_{\b'}^{-1}\ \  ({\text{+ terms of lower degree in
$t_{\Sig_1}$ in case $g=1$)}} \]

We next
restrict $X$ to be a manifold which has the property that the above loops $\bar{a}_i$
on $C\x {\text{\{point\} }}$ bound
vanishing cycles (disks of self-intersection $-1$) in $X\- C$. This means that the boundary of the disk is allowed to move only in $C\x {\text{\{point\} }}$ when computing this self-intersection. For
example, if we let $X$ be the simply connected minimal elliptic surface without multiple fibers and with $\chi = n+1$, $X=E(n+1)$, and let $C$ be a fiber of the genus $n$ fibration on
$E(n+1)$, then the pair $(X,C)$ will satisfy this hypothesis. When this hypothesis is satisfied we will say that $Y'_{1,g,n-1}$ and $(X,C)$ are {\bf{complementary}}. 

Lemma~\ref{rim} shows that the rim tori in $Y'_{1,g,n-1}\- S'$ which come from  
$H_1(S_1\- \{ {\rm{pt}}\};\Z)\to H_1(S';\Z)$ are nullhomologous in $Y'\- S'$. 
Let $R_i$ denote the rim tori in $Z'_{1,g,n-1}(X,C)$ corresponding to the $\bar{a}_i$.  If we assume that  $Y'_{1,g,n-1}$ and $(X,C)$ are complementary, then each $\bar{a}_j$ bounds a disk of self-intersection $-1$ in $X\- C$. Now $a_j$ lies on a fiber $\Sig_{n-1}$ in $Y_{g,n-1}$; so it lies in some  copy of $S^1\x M_{K_{n-1}}$. 
As we have seen in \S2, $a_j$ bounds a punctured torus of self-intersection $+1$ in 
this copy of $S^1\x M_{K_{n-1}}$, and we can clearly make this punctured torus miss a section and a fiber of $S^1\x M_{K_{n-1}}$. In $Z'_{1,g,n-1}(X,C)$ the $(-1)$-disk and the punctured torus glue together to give a torus $U_j$ of self-intersection $0$. The $U_j$ satisfy
\begin{equation*} 
R_i\cdot U_j =\begin{cases} \pm 1,\ \ & j= i \pm 1\\ \ 0 & j \ne i\pm 1\end{cases}
\end{equation*}
Again
as in \S2, if $\k$ is a basic class of $X$ which satisfies $\k\cdot C
= 2n-2$ then we let $\b'_{\k}$ be the class represented by summing
representatives of $\b'$ and $\k$ along $C=S'$, and let
$\ve'_{\k}=\b'_{\k}+2C$.
Then \cite{MST} implies that
\[
\sum_{r_1,\dots, r_n}\sw_{Z'_{1,g,n-1}(X,C)}(\ve'_{\k}+\sum_i r_iR_i)=\sw_{Y'_{1,g,n-1}}(\b')
\cdot\sw_X(\k)= \sw_X(\k) \]
and applying the adjunction inequality to the classes
$U_j$ shows that all $r_i=0$. Thus
for basic classes $\k$ of $X$
which intersect $C$ maximally,
\[ \sw_{Z'_{1,g,n-1}(X,C)}(\ve'_{\k})=\sw_X(\k).\]
Hence $Z'_{1,g,n-1}(X,C)$ has the
same Seiberg-Witten (max) as $Z(X,C,g)$.

\begin{prop} Let $X$ be a symplectic 4-manifold containing
an embedded symplectic surface $C$ of genus $n\ge 2$ and self-intersection
$0$.  Suppose also that $\pi_1(X\- C)=0$ and that $Y'_{1,g,n-1}$ and $(X,C)$ are complementary. Then
$Z'=Z'_{1,g,n-1}(X,C)$ is simply connected, and if the Seiberg-Witten
invariant of $X$ is
\[\sw_X=\sw_{\{X;C,-\}}+\sum_{k\cdot C=2n-2} c_k\,(t_{k}+
(-1)^{\chi(X)}t_{k}^{-1})\]
then the the Seiberg-Witten invariant of
$Z'$ is
\[\sw_{Z'}=\sw_{\{Z';C,-\}}+\sum_{k\cdot C=2n-2}c_k\,(t_{\ve'_{k}}+
(-1)^{\chi(Z')}t_{\ve'_{k}}^{-1})\]
\end{prop}
Whether $Z(X,C,1)$ and $Z'_{1,g,n-1}(X,C)$ are in fact diffeomorphic is a
very interesting question.

More generally, let $L=\{ k_1,\dots, k_n\}$ be a set of positive 
integers, and let 
\[ Y'_{1,g,L}=Y_{1,g}\#_{\Sig_{g,i}=S_{g,i}}\coprod_{i=1}^n Y_{g,k_i} \] 
where $\Sig_{g,i}$ are $n$
genus $g$ fibers of $Y_{1,g}$ and $S_{g,i}$ is the genus $g$ section of
$Y_{g,k_i}$.  Then $Y'_{1,g,L}$ is a symplectic manifold which
contains an embedded symplectic
surface $S'$ of genus $1+\sum k_i$ and self-intersection number $0$ formed from the
sum of a section of $Y_{1,g}$ and fibers of the $Y_{g,k_i}$. The canonical class of
$ Y'_{1,g,L}$ is $\b'=(2\sum k_i)\Sig_g + (2g-2)S'$. 

Now let $X$ be a symplectic 4-manifold with an embedded symplectic
surface $C$ of genus $1+\sum k_i$ and self-intersection $0$. Then we can
form
\[ Z'_{1,g,L}(X,C) = X\#_{C=S'}Y'_{1,g,L} \]

As above, if $X\- C$ is simply connected, then $Z'_{1,g,L}(X,C)$ is also simply
connected. Arguments which are by now familiar show that
 \[ \sw_{Y'_{1,g,n-1}}=t_{\b'}+t_{\b'}^{-1}\ \  ({\text{+ terms of lower degree in
$t_{\Sig_g}$ in case $g=1$)}} \]

In  $Z'_{1,g,L}(X,C)$ we can form classes $\b'_k$ and $\ve'_k$ corresponding to the basic
classes $k$ of $X$ which intersect $C$ maximally. If $\k$ is the
canonical class of $X$ then $\ve'_{\k}$ is the canonical class of
$Z'_{1,g,L}(X,C)$. There is an obvious extension of the definition of `complementarity'  
for $(X,C)$ and $Y'_{1,g,L}$.

\begin{prop} Let $L=\{ k_1,\dots, k_n\}$ be a set of positive 
integers, and let  $X$ be a symplectic 4-manifold containing
an embedded symplectic surface $C$ of genus $1+\sum k_i$ and self-intersection
$0$.  Suppose also that $\pi_1(X\- C)=0$ and that $(X,C)$ and $Y'_{1,g,L}$ are complementary. Then
$Z'= Z'_{1,g,L}(X,C)$ is simply connected, and if the Seiberg-Witten
invariant of $X$ is
\[\sw_X=\sw_{\{X;C,-\}}+\sum_{k\cdot C=2n-2} c_k\,(t_{k}+
(-1)^{\chi(X)}t_{k}^{-1})\]
then the the Seiberg-Witten invariant of
$Z'$ is
\[\sw_{Z'}=\sw_{\{Z';C,-\}}+\sum_{k\cdot C=2n-2}c_k\,(t_{\ve'_{k}}+
(-1)^{\chi(Z')}t_{\ve'_{k}}^{-1})\]
\end{prop}

As we noted above, the hypotheses of this proposition hold for $X$ 
the elliptic surface $E(n)$, $n=2+\sum k_i$, and $C$ a smooth fiber of its 
genus $(n-1)$-fibration.

\section {How to make $Z'_{1,g,L}(X,C)$ nonsymplectic}

Next we show how to modify the manifolds $Z'=Z'_{1,g,L}(X,C)$ in order to
manipulate  the Seiberg-Witten
invariant in such a way that it becomes impossible for the resulting
manifold to admit a symplectic structure.  We first
identify families of nullhomologous tori in $Z'$ with self-intersection
number $0$  upon which we will perform surgery.
The manifold
$Y_{1,g}= S^1\x M_{K_g}$ is the total space of a fiber bundle $p: Y_{1,g}\to B$
with a genus $1$ base and genus $g$ fiber, $\Sig_g$. Let $S_1^1$ and $S_2^1$ be embedded circles in the base which generate $\pi_1$ and such that the $p^{-1}(S_1^1)\cong S^1_1\x \Sig_g$ and $p^{-1}(S^1_2)\cong M_{K_g}$. We may suppose that $S_1^1$ and $S_2^1$ intersect in a single point, $x_0$. For any closed embedded loop $\a$ in $p^{-1}(x_0)$, there is a torus
$\L(\a)=S^1_1\x \a\C Y_{1,g}$.

For the fiber sum $Y'_{1,g,L}=Y_{1,g}\#_{\Sig_{g,i}=S_{g,i}}\coprod_{i=1}^n Y_{g,k_i}$,  let $b_i$, $i=1,\dots, n$ be the points in $B$ whose fibers $\Sig_{g,b_i}$ are identified with the sections $S_{g,i}$ of the $Y_{g,k_i}$.  We may assume that $b_i\in B\- (S^1_1\cup S^1_2)$ for each $i$, and so $\L(\a)$ represents a class in $H_2(Y'_{1,g,L};\Z)$. This group is generated by $\Sig_g$, $S'$, the rim tori and vanishing classes in the $Y_{g,k_i}$, and the rim tori and vanishing classes which arise from the $\Sig_{g,b_i}$. However, since the $\Sig_{g,b_i}$ are identified with sections $S_{g,i}$ in $Y'_{1,g,L}$, these last rim tori  are nullhomologous in $Y'_{1,g,L}$ and there are no corresponding vanishing classes.
It is clear that $\L(\a)$ is disjoint from all the rest of these surfaces; so this means that in $H_2(Y'_{1,g,L};\Z)$, $\L(\a)$ is nullhomologous. We shall be interested in the case where 
$\a=a_i$, the loops giving the basis for $H_1(\Sig_g;\Z)$ which is described in \S2. 

We now wish to perform surgeries on collections of the $\L(a_i)$. In $S^3\- K_g$, a positive push-off $a_i^+$ of $a_i$ bounds a punctured torus in the complement of $K_g$. Since  $p^{-1}(S^1_2)$ is diffeomorphic to $M_{K_g} = (S^3\- (K_g\x D^2)) \cup (S^1\x D^2)$, we see that there is a push-off of $a_i$ onto $\bd(a_i\x D^2)$ which bounds a punctured torus in $M_{K_g}\- (a_i\x D^2)$. Denote this push-off by $\ha_i$. It defines a `$0$-framing' for $a_i$, and we use it to express the boundary of a tubular neighborhood of $\L(a_i)$ as
$S^1_1\x \ha_i\x \bd D^2$. In general, remove this tubular neighborhood
from $Z$ and reglue it so that (in homology) $\bd D^2$ is sent to
$m\, \ha_i+\ell \,\bd D^2$ and $S^1_1$ to $S^1_1$ on $\bd(Z'\- (S^1_1\x \ha_i \x
D^2))$. This gives $\ell/m$-surgery
on $\L(\a_i)$. We are interested in $(-1/m)$-surgery for all $m\ne 0$.

For $m\ne0$, $(-1/m)$-surgery on $a_1$ in Figure~2
turns $K_g$ into a different knot
$K_g(m)$ in $S^3$. The construction gives an obvious genus
$g$ Seifert surface $\Sig$ for $K_g(m)$. (See Figure 4.) It is then
clear that $(-1/m)$-surgery on $\L(a_1)$ in  $S^1\x  M_{K_g}$ gives $S^1\x M_{K_g(m)}$. 
Performing this surgery on  $\L(a_1)\C Y'_{1,g,L}$ gives us a new manifold 
\[ Y'_{1,g,L}(m)= S^1\x M_{K_g(m)}\,\#_{\Sig_{g,i}=S_{g,i}}\coprod_{i=1}^n Y_{g,k_i} \]
and
\[ Z'(m) = Z'_{1,g,L}(X,C;m) =X\#_{C=S'}Y'_{1,g,L}(m) \]
since $S'$ is still defined using the torus $S^1\x \{ {\text{meridian}}\}$ in $S^1\x M_{K_g(m)}$
and the fibers $\Sig_{k_i}$ in the $Y_{g,k_i}$. The argument of Propostion~\ref{pi1Z'} shows that if $X\- C$ is simply connected, then so is $Z'(m)$. 

We next wish to calculate the Seiberg-Witten invariant of $Z'(m)$. 
Using the Mayer-Vietoris sequence for
$Z' = (Z' \- (\L(a_1)\x D^2)) \cup (\L(a_1)\x D^2)$ and similarly for $Z'(m)$, we see that we may identify the homology groups of these manifolds. Let $\L(a_1;m)$ be the core torus of $\L(a_1)\x D^2$ in $Z'(m)$. Since we are performing $(-1/m)$-surgery on the nullhomologous torus $\L(a_1)$, the torus $\L(a_1;m)$ is also nullhomologous. Let $\hZ'$ denote the
result of the surgery on $\L(a_1)$ for which $\bd D^2$ is sent to
$\ha_1$ and $S^1_1$ to $S^1_1$ on $\bd(Z'\- (S^1_1\x \ha_1 \x
D^2))$.  (This corresponds
to $(0/1)$-surgery on $\L(a_1)$.) The intersection number of a
basic class with a self-intersection $0$ torus must be $0$; so the basic
classes of $Z'$ may also be viewed as homology classes of $Z'(m)$ which
have trivial intersection with the image torus $\L(a_1;m)$.
In $\hZ'$ there is the additional class
$\t$ which is represented by the torus $S^1_1\x \{\text{point in $a_1$}\}\x\bd D^2$  on
$\bd(Z'\- (\L(a_1)\x D^2))=\bd(Z'\- (S^1_1\x a_1\x D^2))$.
Since $\L(a_1)$ is nullhomologous, the surgery formula of Morgan, Mrowka, Szabo \cite{MMS} and Taubes \cite{T3} (see also \cite{KL4M}) states
\[ \sw_{Z'(m)}(\z)=\sw_{Z'}(\z)-m\sum_{i,j}\sw_{\hZ'}(\z+j\t). \]

For each basic class $k$ of $X$ which satisfies $k\cdot C=2\sum k_i$ we have classes $\b'_k$ and $\ve'_k$ in $Z'$. There are also such classes in $Z'(m)$, and to avoid further notational complexities, we shall still denote them by $\b'_k$ and $\ve'_k$.

\centerline{\unitlength 1cm
\begin{picture}(12,11)
\put (2,10.5){\line(-1,0){1}}
\put (1,10.14){\oval(1.2,.75)[tl]}
\put (2,10.5){\line(3,-4){1.2}}
\put (3,10.5){\line(-3,-4){.4}}
\put (1.8,8.89){\line(-1,0){1.415}}
\put (.385,10.15){\line(0,-1){1.25}}
\put (.385,7.14){\framebox(1.415,1.75){}}
\put (.8, 8.5){\Small{$m$}}
\put (.7, 8.1){\Small{full}}
\put (.71,7.7){\Small{RH}}
\put (.55, 7.3){\Small{twists}}
\put (3.2,8.89){\line(0,-1){1.75}}
\put (2.4,9.69){\line(-3,-4){.6}}
\put (1.8,7.14){\line(3,-4){1.4}}
\put (3.2,7.14){\line(-3,-4){.6}}
\put (1.8,5.27){\line(3,4){.6}}
\put (2.4,4.7){$\bullet$}
\put (2.4,4.03){$\bullet$}
\put (2.4,3.36){$\bullet$}
\put (1.8,3){\line(3,-4){1.4}}
\put (3.2,3){\line(-3,-4){.6}}
\put (1.8,1.13){\line(3,4){.6}}
\put (1.8,1.13){\line(-1,0){1}}
\put (1,1.52){\oval(1.2,.75)[bl]}
\put (.385,7.14){\line(0,-1){5.65}}
\put (3.2,1.13){\line(1,0){1}}
\put (4,1.52){\oval(1.2,.75)[br]}
\put (3,10.5){\line(1,0){1}}
\put (4,10.14){\oval(1.2,.75)[tr]}
\put (4.59,10.15){\line(0,-1){8.67}}
\put (1.6,.5){Figure 4}
\put (8,6){\oval(3,4)[l]}
\put (9,6){\oval(3,4)[r]}
\put (8,8){\line(3,-4){.4}}
\put (9,6.66){\line(-3,4){.4}}
\put (9,6.66){\line(-3,-4){.4}}
\put (8,6.66){\line(3,4){1}}
\put (8,6.66){\line(3,-4){1}}
\put (8,5.33){\line(3,4){.4}}
\put (8,5.33){\line(3,-4){1}}
\put (9,5.33){\line(-3,-4){.4}}
\put (8,4){\line(3,4){.4}}
\put (7.85,3.5){Figure 5}
\end{picture}}

The knot $K_g(-1)$, obtained from $(+1)$-surgery on $a_1$, has a Seifert surface
of genus $g-1$, one less than that of $K_g$. Figure 5 exhibits
the case $g=1$. In $M_{K_g(-1)}$ and therefore in  $Y'_{1,g,L}(-1)$ there is a corresponding closed surface $\G_{g-1}$. Calculating as in \S2 we find that the only possible basic classes of $Y'_{1,g,L}(-1)$ are $\pm((2\sum k_i)\Sig_g +(2g-2)S')$ (and lesser multiples of $\Sig_g$ in case $g=1$). However, $S'\cdot \G_{g-1}=1$ and $\Sig_g\cdot\G_{g-1}=0$. In case $g>1$, the adjunction inequality, applied to $\G_{g-1}$, leads directly to a contradiction if we assume that 
$\pm((2\sum k_i)\Sig_g +(2g-2)S')$ is a basic class. In case $g=1$, $\G_0$ is a sphere of self-intersection $0$, and it is essential since $S'\cdot\G_0=1$ (\cite{Turkey}). Thus the Seiberg-Witten invariant of $Y'_{1,g,L}(-1)$ vanishes in this case as well. Using \cite{MST}, we get $\sw_{Z'(-1)}=0$, and
applying the above formula we obtain
\[ 0=\sw_{Z'}(\z)+\sum_i\sw_{\hZ'}(\z+i\t) \]
for each basic class $\z$ of $Z'$, and hence:

\begin{thm} \label{b} Let $L=\{ k_1,\dots, k_n\}$ be a set of positive 
integers, and let  $X$ be a symplectic 4-manifold containing
an embedded symplectic surface $C$ of genus $1+\sum k_i$ and self-intersection
$0$.  Suppose also that $\pi_1(X\- C)=0$ and that $(X,C)$ and $Y'_{1,g,L}$ are complementary. Then for any positive integer $m$, $Z'(m)=Z'_{1,g,L}(X,C;m)$ is simply connected, and if the Seiberg-Witten invariant of $X$ is
\[ \sw_X=\sw_{\{X;C,-\}}+\sum\limits_{k\cdot C=2n-2} c_k\,(t_{k}+
(-1)^{\chi(X)}t_{k}^{-1})\]
then the the Seiberg-Witten invariant of
$Z'(m)$ is
\[\sw_{Z'(m)}=\sw_{\{Z'(m);C,-\}}+(m+1)\sum_{k\cdot C=2n-2} c_k\,(t_{\ve'_{k}}+
(-1)^{\chi(Z'(m))}t_{\ve'_{k}}^{-1}).\]
\end{thm}
\noindent Of course, it follows that  if $m\ne 0$, $Z'(m)$ can admit no
symplectic structure.

It is clear that nothing special is gained by working with $a_1$ in this construction. If 
$\bm =(m_1,\dots,m_{2g})$ where the $m_i$ are nonnegative integers, and
we perform $(-1/{m_i})$-surgeries on (push-offs) of the $a_i$ in $M_{K_g}$ we obtain $0$-surgery on a new knot $K_g(\bm)$. Our construction then leads us to manifolds
$Z'(\bm) = Z'_{1,g,L}(X,C;\bm)$ As long as $X\- C$ is simply connected, 
$Z'(\bm)$ will also be simply connected, and, with notation as above,

\begin{cor}
$ \sw_{\{Z'(\bm);C,{\text\em{max}\} }}=(\prod\limits_{i=1}^{2g}(m_i+1))\cdot \sum\limits_{k\cdot C=2n-2} c_k\,(t_{\ve_{k}}+
(-1)^{\chi(Z'(\bm))}t_{\ve_{k}}^{-1})$.
\end{cor}

In the special case where $X$ is the result of knot surgery on the
K3-surface with a knot of genus $\sum k_i$ and $C$ is the surface described in the paragraph below Proposition~\ref{Z}, $(X,C)$ and $Y'_{1,g,L}$ are complementary, and we get
\[ \sw_ {Z'(\bm)} = (\prod\limits_{i=1}^{2g}(m_i+1))\cdot (t_F^{n-1}\pm t_F^{1-n}).\]
So, for example, for any of the many
choices of $\bm$ with corresponding products
$\prod(m_i+1)$ equal, we have manifolds which cannot be
distinguished by means of their Seiberg-Witten invariants.  One can also vary 
$L=\{ k_1,\dots,k_n \}$.  

For the record, we note following facts about the characteristic
numbers $c_1^2(Z'(\bm))$ and
$\chi(Z'(\bm))=\frac{1}{4}(\sign(Z'(\bm))+e(Z'(\bm)))$.

\begin{prop} The characteristic numbers of $Z'(\bm)=Z'_{1,g,L}(X,C;\bm)$ are:
\begin{eqnarray*}
c_1^2(Z'(\bm)=c_1^2(X)+8(n(g-1)+\sum k_i),\hspace{.1in}{\text{and}}\\
\chi(Z'(\bm))=\chi(X)+(n(g-1)+\sum k_i). \end{eqnarray*}
 \end{prop}

\section {Symplectic fibrations}

A theorem of Parshin and Arakelov \cite{Par,Ar} (see also \cite{JY})
states that given fixed
$g\ge2$ and a finite set of points $S\C\PO$, there are at most
finitely many holomorphic
fibrations over $\PO$ whose generic fiber  is a Riemann surface of
genus $g$ and whose
singular fibers have image in $S$.   Interest in fibrations of
symplectic manifolds has been
rekindled by the theorem of Simon Donaldson \cite{SKD} which states
that (after blowing up)
each symplectic
$4$-manifold admits a locally holomorphic Lefschetz fibration over
$S^2$. In this section we shall give examples which show that the
Parshin-Arakelov  Theorem
has no analogue in the symplectic category. (See Theorem~\ref{LF}
below.) The authors have
described more complicated examples exhibiting this phenomenon in
\cite{surfaces}.

The Lefschetz fibrations which we have in mind live on the homotopy
elliptic surfaces of
\cite{KL4M}. These manifolds,
$E(n)_K$, are built from knot surgery on the simply connected,
minimally elliptic surface
$E(n)$ without multiple fibers and of holomorphic Euler
characteristic $n$. If $K$ is a
nontrivial fibered knot, then
$E(n)_K$ is a symplectic $4$-manifold which admits no complex
structure (nor does
$E(n)_K$ with the opposite orientation).

The elliptic surface $E(n)$ is the double branched cover of $\SS$
with branch set equal to
four disjoint copies of $S^2\x \{{\rm{pt}}\}$ together with
$2n$ disjoint copies of $\{{\rm{pt}}\}\x S^2$. The resultant branched
cover has $8n$ singular
points (corresponding to the double points in the branch set), whose
neighborhoods are cones
on
${{\mathbf{RP}}^{\,3}}$. These are desingularized in the usual way,
replacing their
neighborhoods with cotangent bundles of $S^2$. The result is $E(n)$.
The horizontal and
vertical fibrations of
$\SS$ pull back to give fibrations of
$E(n)$ over $\PO$. A generic fiber of the vertical fibration is the
double cover of $S^2$,
branched over $4$ points --- a torus. This describes an elliptic fibration of
$E(n)$. The generic fiber of the horizontal fibration is the double cover of
$S^2$, branched over $2n$ points, and this gives a genus $n-1$ fibration on
$E(n)$. This genus $n-1$ fibration has four singular fibers which are
the preimages of the
four $S^2\x \{{\rm{pt}}\}$'s in the branch set  together with the
spheres of self-intersection
$-2$ arising from desingularization.  The generic fiber
$T$ of the elliptic fibration meets a generic fiber
$\Sig_{n-1}$ of the horizontal fibration in two points,
$\Sigma_{n-1}\cdot T=2$.

Let $K$ be a fibered knot of genus $g$, and fix a generic elliptic
fiber $T_0$ of
$E(n)$. Then in the knot surgery manifold
\[ E(n)_K = (E(n)\- (T_0\x D^2)) \cup (S^1\x (S^3\- N(K)) ,\] each
normal $2$-disk to $T_0$
is replaced by a fiber of the fibration of $S^3\- N(K)$ over
$S^1$. Since
$T_0$ intersects each generic horizontal fiber twice, we obtain a
`horizontal' fibration
\[ h: E(n)_K\to\PO\]  of genus $2g+n-1$.

This fibration also has four singular fibers arising from the four
copies of $S^2\x
\{{\rm{pt}}\}$ in the branch set of the double cover of $\SS$. Each
of these gets blown up at
$2n$ points in $E(n)$, and the singular fibers each  consist of a genus
$g$ surface $\Sig_g$ of self-intersection $-n$ and multiplicity $2$
with $2n$ disjoint
$2$-spheres of self-intersection $-2$, each meeting $\Sig_g$
transversely in one point. The
monodromy around each singular fiber is (conjugate to) the diffeomorphism of
$\Sig_{2g+n-1}$ which is the deck transformation $\eta$ of the double cover of
$\Sig_g$, branched over $2n$ points. Another way to describe $\eta$
is to take the
hyperelliptic involution $\o$ of $\Sig_{n-1}$ and to connect sum copies of
$\Sig_g$ at the two points of a nontrivial orbit of $\o$. Then $\o$
extends to the involution
$\eta$ of $\Sig_{2g+n-1}$.

The fibration which we have described is not Lefschetz since the
singularities are not simple
nodes. However, it can be perturbed locally to be Lefschetz:

\begin{lem} \label{Mng} Any symplectic fibration on a $4$-manifold
with singular fibers
equivalent to those of $h: E(n)_K\to\PO$ can be locally deformed to a
Lefschetz fibration.
\end{lem}
\begin{proof} It suffices to find a holomorphic model for the
singular fibers of
$h$, since these can be locally deformed to a Lefschetz fibration by
complex Morse theory.
The model is built from a branched double cover of $\Sig_g\x S^2$.
This time the branch set
consists of two disjoint copies of $\Sig_g\x \{{\rm{pt}}\}$ and $2n$
disjoint copies of
$\{{\rm{pt}}\}\x S^2$. After desingularizing as above, one obtains a
complex surface
$M(n,g)$ with a holomorphic (horizontal) fibration of genus
$2g+n-1$ and with a pair of singular fibers exactly of the type of
the singular fibers of $h$.
\end{proof}

Returning to our examples, we have:

\begin{thm}{\label{LF}} If $K$ is a fibered knot whose fiber has genus
$g$, then $E(n)_K$  admits a locally holomorphic fibration (over
$\PO$) of genus $2g+n-1$ which has exactly four singular fibers.
Furthermore, this fibration
can be deformed locally to be Lefschetz.
\end{thm}

There are qualitative differences in the fibrations on $E(n)_K$
between $n=1$ and $n>1$. From
their local description in Lemma~\ref{Mng}, $E(1)_K$ contains
singular fibers which are
reducible; i.e. vanishing cycles obtained from Dehn twists about
separating curves. This is
due to the fact that it is built from a genus $0$ fibration on
$E(1)$.  Siebert and Tian
\cite{ST} have shown that any genus~2 Lefschetz fibration with only
irreducible singular
fibers must be holomorphic. Using this, Auroux \cite{AR} has shown
that any genus~2
Lefschetz fibration is stably holomorphic; i.e. the fiber sum of any
genus~2 Lefschetz
fibration with sufficiently many copies of the rational genus~2
Lefschetz fibration $X_2$ with
20 irreducible singular fibers is isomorphic to a holomorphic
fibration. Thus, when $K$ is
either the trefoil or figure $8$ knot (genus 1), then
$E(1)_K$ admits a genus~$2$ Lefschetz fibration and $E(1)_K$ admits
no complex structure.
However, the fiber sum of $E(1)_K$ with four copies of $X_2$ is
diffeomorphic to a complex
surface.  The case of $E(n)_K$, $n>1$, is different. From their local
description in
Lemma~\ref{Mng}, the fibrations in Theorem~\ref{LF} are all of genus
larger than 2 and have
only irreducible singular fibers:

\begin{prop} For $n\ge 2$, the genus $2g+n-1$ Lefschetz fibrations, described above, on the manifolds $E(n)_K$ have no reducible singular fibers.\end{prop}
\begin{proof}  Our argument determines the vanishing cycles of the Lefschetz fibration, and shows that they must all be nonseparating, or equivalently the corresponding singular fibers must be irreducible. It is easy to see that for any genus $G$ Lefschetz fibration on a manifold $X$ over $S^2$, the euler number of $X$ is given by $e(X)=s-4G+4$, where $s$ is the number of singular fibers. Since $e(E(n)_K) =12n$, our Lefschetz fibration has $16n+8g-8$ singular fibers. This fibration is obtained by a local deformation of a fibration $\xi$ on $E(n)_K$ which has $4$ singular fibers. Each of these $4$ singular fibers contributes $4n+2g-2$ singular fibers and therefore $4n+2g-2$ vanishing cycles to the Lefschetz fibration.

 By construction, the Lefschetz fibration on $E(n)_K$ is obtained from the standard genus $n-1$ hyperelliptic Lefschetz fibration on $E(n)$. This fibration has $16n-8$ singular fibers, and, for $n\ge 2$, it is well-known (see, e.g. \cite{GS}) that they are all irreducible. 

We next take account of the fact that the singular fibers of $\xi$ have genus $g$, whereas the nonsingular fibers have genus $2g+n-1$. Vanishing cycles must account for the corresponding reduction of the first betti numbers of the fibers. Already we have taken into account the reduction in genus from $n-1$ to $0$, since this occurs in the hyperelliptic fibration. It remains to find $2g$ more vanishing cycles for each singular fiber of $\xi$. This accounts for the remaining $8g$ vanishing cycles, and they all must be nonseparating, since a separating vanishing cycle does not reduce the first homology of the fiber.
\end{proof}

 None of these examples, for $n\ge 2$, are
hyperelliptic. Thus,
while it is a conjecture of Siebert and Tian \cite{ST} that a
hyperelliptic Lefschetz
fibration with no reducible fibers is indeed holomorphic, the
examples constructed in
Theorem~\ref{LF} show that hyperellipticity is a necessary hypothesis.

There is another way
to view these constructions. Let $M(n,g)$ be the complex surface
arising in the proof of
Lemma~\ref{Mng}. Once again, this  manifold carries a pair of
fibrations. There is a genus
$2g+n-1$ fibration over $S^2$  and an $S^2$ fibration over $\Sig_g$.

Consider first the $S^2$ fibration. This has $2n$ singular fibers,
each of which consists of
a smooth 2-sphere $E_i$, $i=1,\dots, 2n$, of  self-intersection $-1$
and multiplicity $2$,
together with a pair of disjoint spheres of self-intersection
$-2$, each intersecting $E_i$ once transversely. If we blow down
$E_i$ we obtain again an $S^2$ fibration over $\Sig_g$, but the $i$th
singular fiber  now
consists of a pair of 2-spheres of self-intersection $-1$ meeting
once, transversely. Blowing
down one of these gives another $S^2$ fibration over $\Sig_g$, with
one less singular fiber.
Thus blowing down $M(n,g)$ $4n$ times results in  a manifold which is
an $S^2$ bundle over
$\Sig_g$. This means that (if $n>0$) $M(n,g)$ is diffeomorphic to
$(S^2\x \Sig_g)\# 4n\,\CPb$.

The genus $2g+n-1$ fibration on $M(n,g)$ has 2 singular fibers. As
above, these fibers consist
of a genus $g$ surface $\Sig_g$ of self-intersection $-n$ and
multiplicity $2$ with $2n$
disjoint $2$-spheres of self-intersection $-2$, each meeting $\Sig_g$
transversely in one
point. The monodromy of the  fibration around each of these fibers is
the deck transformation
of the double branched cover of
$\Sig_g$. This is just the map $\eta$ described above.

Let $\vp$ be a diffeomorphism of $\Sig_g\- D^2$ which is the identity
on the boundary. For
instance, $\vp$ could be the monodromy of a fibered knot of genus
$g$. There is an induced diffeomorphism $\Phi$ of \,
      $\Sig_{2g+n-1}=\Sig_g\#\Sig_{n-1}\#\Sig_g$ which is given by $\vp$
on the first
$\Sig_g$ summand and by the identity on the other summands. Consider
the twisted fiber sum
\[ M(n,g)\#_{\Phi}M(n,g) = \{M(n,g)\- (D^2\x \Sig_{2g+n-1})\}
\cup_{{\text{id}}\x\Phi} \{M(n,g)\- (D^2\x \Sig_{2g+n-1})\} \] where
fibered neighborhoods of
generic fibers $\Sig_{2g+n-1}$ have been removed from the two copies
of $M(n,g)$, and they
have been glued by the diffeomorphism
${{\text{id}}\x\Phi}$ of $S^1\x \Sig_{2g+n-1}$.

In the case that  $\vp$ is the monodromy of a fibered knot $K$, we claim that
$M(n,g)\#_{\Phi}M(n,g)$ is the manifold $E(n)_K$ with the genus
$2g+n-1$ fibration described above. To see this, we view $S^2$ as the
base of the horizontal
fibration. Then it suffices to check that the total monodromy map
$\pi_1(S^2\- 4\, {\text{points}})\to {\text{Diff}}(\Sig_{2g+n-1})$ is
the same for  each. It
is not difficult to see that if we write the generators of
$\pi_1(S^2\- 4\, {\text{points}})$ as $\a$, $\b$, $\g$ with $\a$ and
$\b$ representing loops
around the singular points of, say, the image of  the first copy of
$M(n,g)$ and basepoint in
this image, and $\g$ a loop around a  singular point in the image of
the second $M(n,g)$ then
the monodromy map $\mu$ satisfies
$\mu(\a)=\eta$, $\mu(\b)=\eta$ and $\mu(\g)$ is
$\vp\oplus\o\oplus\vp^{-1}$, expressed as a
diffeomorphism of $\Sig_g\#\Sig_{n-1}\#\Sig_g$. That  this is also
the monodromy of $E(n)_K$
follows directly from its construction.

\section {Construction III}

Let $K_1$, $K_2$ be fibered knots of genus $g$, and fix $n\ge 1$.
Then as in the previous
section there are genus $2g+n-1$ fibrations on the manifolds
$E(n)_{K_i}$. The manifolds $E(n)_{K_i}$ do not admit torus
fibrations, but it makes sense to
speak of `elliptic fibers' $T_i$ in $E(n)_{K_i}$, meaning those remaining from
$E(n)$ after knot surgery. As above, $T_i$ intersects the horizontal fiber
$\Sig_{2g+n-1}$ in two positive intersection points.

Let $Y=Y(n;K_1,K_2)$ denote the fiber sum
$E(n)_{K_1}\,\#_{\Sig_{2g+n-1}}E(n)_{K_2}$. Then $Y$ is a symplectic
manifold with
$c_1^2(Y)=16g+8n-16$. Furthermore, $Y$ is simply connected because in
$E(n)_{K_i}$ a remnant of a singular fiber of the elliptic fibration
on $E(n)$ provides a
2-sphere (of self-intersection $-2$) which intersects the fiber
$\Sig_{2g+n-1}$ in one point. Tori
$T_i\C E(n)_{K_i}$, $i=1,2$ which intersect $\Sig_{2g+n-1}\x S^1$ in
identified circles will
glue together to form a genus 3 surface in $Y$.   Let
$\t$ denote the homology class in $Y$ which is represented by this
surface. Hence
$\t\cdot\Sig_{2g+n-1}=2$. The canonical class of $Y$ is
$K_Y=(2g+n-2)\,\t+2\,\Sig_{2g+n-1}$, and any basic class $\k$ of $Y$ satisfies
$\k^2=16g+8n-16$.

Straightforward adjunction inequality arguments show that the only
basic classes of
$Y$ are those of form $\pm K_Y +R$ where $R$ lies in the subgroup of
$H_2(Y;\Z)$ generated by rim tori of $\Sig_{2g+n-1}$. Let $\o$ denote
a  symplectic form on $Y$
obtained as a result of symplectic fiber sum. Then all rim tori of
$\Sig_{2g+n-1}$ are Lagrangian with respect to $\o$. Invoking
\cite{T}, we see that
$\int_{K_Y}\o$ is the unique maximal value of $\int_{\k}\o$ among all
basic classes $\k$  of
$Y$. But because rim tori are Lagrangian,
\[ \int_{K_Y+R}\o = \int_{K_Y}\o \] This means that $K_Y + R$ is not basic if
$R\ne 0$. Thus the only basic classes of
$Y$ are $\pm K_Y$, and $\sw_Y=t_K+(-1)^n\,t_K^{-1}$.

If $\{K_i\}$ is a family of genus $g$ fibered knots, then the above
discussion shows that
manifolds $Y(n; K_i, K_j)$ cannot be distinguished on the  basis of
their Seiberg-Witten
invariants.

\end{document}